\documentclass[11pt, a4paper]{article}
\usepackage{amsmath, amssymb}
\usepackage{amsthm}
\usepackage{hyperref}

\newtheorem{lemma}{Lemma}[section]

\newtheorem{theorem}[lemma]{Theorem}
\theoremstyle{definition}

\newtheorem{remark}[lemma]{Remark}
\newcommand{\aire}{\mathrm{Area}(M^2,g)}
\newcommand{\R}{\mathbb{R}}
\newcommand{\B}{\mathcal{B}}
\newcommand{\SSS}{\mathbb{S}}
\newcommand{\C}{\mathbb{C}}

\newcommand{\Spinc}{\mathrm{Spin^c}}

\newcommand{\bq}{\begin{equation}}
\newcommand{\bqw}{\begin{equation*}}
\newcommand{\eq}{\end{equation}}
\newcommand{\eqw}{\end{equation*}}
\newcommand{\vo}{\mathrm{v}}
\newcommand{\ov}{\overline}

\newcommand{\vol}{\mathrm{vol}}

\renewcommand{\i}{\mathrm{i}}

\renewcommand{\Re}{\mathrm{Re}}
\newcommand{\Spin}{\mathrm{Spin}}

\newcommand{\cercle}{\mathbb{S}}
\newcommand{\End}{\mathrm{End}}
\newcommand{\iid}{\mathrm{Id}\,}
\begin{document}
\title{{\bfseries The Hijazi inequalities on  complete Riemannian $\Spinc$ manifolds}}

\author{Roger NAKAD}




\maketitle
\begin{center}
Institut \'Elie Cartan, Universit\'e Henri Poincar\'e, Nancy I, B.P 239\\
54506 Vand\oe uvre-L\`es-Nancy Cedex, France.
\end{center}
\begin{center}
 {\bf nakad@iecn.u-nancy.fr}
\end{center}
\vskip 0.5cm
\begin{center}
 {\bf Abstract}
\end{center}
In this paper, we extend the Hijazi type inequality, involving the Energy-Momentum tensor, to the eigenvalues of the Dirac operator on  complete Riemannian $\Spinc$ manifolds without boundary and of finite volume. Under some additional assumptions, using the refined Kato inequality, we prove the Hijazi type inequality for elements of the essential spectrum. The limiting cases are also studied.
\section{Introduction}
On a compact Riemannian $\Spinc$ manifold $(M^n, g)$ of dimension $n\geqslant2$, any eigenvalue $\lambda$ of the Dirac operator satisfies the Friedrich type inequality \cite{8, friedrich}:
\begin{eqnarray}
 \lambda^2 \geqslant \frac{n}{4(n-1)}\inf_M (S -c_n\vert\Omega\vert),
\label{spinfr}
\end{eqnarray}
where $S$ denotes the scalar curvature of $M$, $c_n =2 [\frac n2]^\frac 12$ and $i\Omega$ is the curvature form of the connection on the line bundle given by the $\Spinc$ structure. Equality holds if and only if the eigenspinor $\psi$ associated with the first eigenvalue $\lambda_1$ is a $\Spinc$ Killing spinor, i.e., for every $X\in \Gamma(TM)$ the eigenspinor $\psi$ satisfies
\begin{eqnarray}
\left\{
\begin{array}{l}
\nabla_X\psi = -\frac{\lambda_1}{n}X\cdot\psi,\\
\Omega\cdot\psi =i\frac {c_n}{2}\vert\Omega\vert \psi.
\end{array}
\label{lamba}
\right.
\end{eqnarray}
Here $X\cdot\psi$ denotes the Clifford multiplication and $\nabla$ the spinorial Levi-Civita connection \cite{5, 16}. In \cite{article1}, it is shown that on a compact Riemannian $\Spinc$ manifold any eigenvalue $\lambda$ of the Dirac operator to which is attached an eigenspinor $\psi$ satisfies the Hijazi type inequality \cite{12} involving the Energy-Momentum tensor and the scalar curvature:
\begin{eqnarray}
 \lambda^2 \geqslant \inf_M (\frac 14 S -\frac{c_n}{4} \vert\Omega\vert + \vert \ell^\psi\vert^2),
\label{oussamascal}
\label{oussamascal11}
\end{eqnarray}
where $\ell^\psi$ is the field of symmetric endomorphisms associated with the field of quadratic forms denoted  by $T^\psi$, called the Energy-Momentum tensor. It is defined on the complement set of zeroes of the eigenspinor $\psi$, for any vector field $X$ by 
$$T^\psi (X)= \Re < X\cdot\nabla_X \psi,\frac{\psi}{\vert\psi\vert^2}>.$$
Equality holds in (\ref{oussamascal}) if and only, for all $X \in \Gamma(TM)$, we have
\begin{eqnarray}
\left\{
\begin{array}{l}
\nabla_X\psi = -\ell^\psi(X)\cdot\psi,\\
\Omega\cdot\psi = i\frac{c_n}{2} \vert\Omega\vert \psi,
\end{array}
\right.
\label{casee}
\label{morel1}
\end{eqnarray}
where $\psi$ is an eigenspinor associated with the first eigenvalue $\lambda_1$. By definition, the trace $tr(\ell^\psi)$ of $\ell^\psi$, where $\psi$ is an eigenspinor associated with an eigenvalue $\lambda$, is equal to $\lambda$. Hence, Inequality (\ref{oussamascal}) improves Inequality (\ref{spinfr}) since by the Cauchy-Schwarz inequality, $\vert \ell^\psi\vert^2 \geqslant \frac{(tr(\ell^\psi))^2}{n} = \frac{\lambda^2}{n}$. It is also shown that the sphere equipped with a special $\Spinc$ structure satisfies the equality case in (\ref{oussamascal}) but equality in (\ref{spinfr}) cannot occur.\\ \\
In the same spirit as in \cite{13}, A. Moroianu and M. Herzlich (see \cite{8}) generalized the Hijazi inequality \cite{13}, involving the first eigenvalue of the Yamabe operator $L$, to the case of compact $\Spinc$ manifolds of dimension $n\geqslant 3$: Any eigenvalue $\lambda$ of the Dirac operator satisfies
\begin{eqnarray}\label{herr}
\lambda^2 \geqslant \frac{n}{4(n-1)}\mu_1,
\end{eqnarray}
where $\mu_1$ is the first eigenvalue of the perturbed Yamabe
operator defined by $L^\Omega = L- c_n\vert\Omega\vert_g = 4\frac{n-1}{n-2}\bigtriangleup + S  -c_n\vert\Omega\vert_g.$ The limiting case of (\ref{herr}) is  equivalent to the limiting case in (\ref{spinfr}). The Hijazi inequality \cite{12}, involving the Energy-Momentum tensor and the first eigenvalue of the Yamabe operator, is then proved by the author in \cite{article1} for compact $\Spinc$ manifolds. In fact, any eigenvalue of the Dirac operator to which is attached an eigenspinor $\psi$ satisfies 
\begin{eqnarray}
\lambda ^2 \geqslant \frac 14 \mu_1 +\inf_M\ \vert \ell^{\psi}\vert^2.
\label{fin}
\label{ega}
\end{eqnarray}
Equality in (\ref{ega}) holds if and only, for all $X \in \Gamma(TM)$, we have 
\begin{eqnarray}
\left\{
\begin{array}{l}
\overline\nabla_X\overline\varphi = -\ell^{\overline\varphi}(X)\ \overline\cdot\ \ \overline\varphi,\\
\Omega\cdot\psi = i\frac{c_n}{2}\vert\Omega\vert_g\psi,
\end{array}
\right.
\label{case}
\end{eqnarray}
 where $\overline\varphi = e^{-\frac{n-1}{2}u}\overline\psi$, the spinor field $\overline \psi$ is the image of $\psi$ under the isometry between the spinor bundles of $(M^n, g)$ and $(M^n, \overline g = e^{2u}g)$ and $\psi$ is an eigenspinor associated with the first eigenvalue $\lambda_1$ of the Dirac operator. Again, Inequality (\ref{fin}) improves Inequality (\ref{herr}). In this paper we examine these lower bounds on open manifolds, and especially on complete Riemannian $\Spinc$ manifolds. We prove the following:
\begin{theorem}
Let $(M^n,g)$ be a complete Riemannian $\Spinc$ manifold of finite volume. Then any eigenvalue $\lambda$ of the Dirac operator to which is attached an eigenspinor $\psi$ satisfies the Hijazi type inequality (\ref{oussamascal11}). 
Equality holds if and only if the eigenspinor associated with the first eigenvalue $\lambda_1$ satisfies (\ref{morel1}).
\label{ththth}
\end{theorem}
The Friedrich type inequality (\ref{spinfr}) is derived for complete Riemannian $\Spinc$ manifolds of finite volume and equality also holds if and only if the eigenspinor associated with the first eigenvalue $\lambda_1$ is a Killing $\Spinc$ spinor.
This was proved by N. Grosse in \cite{nadine} and \cite{nadine1} for complete spin manifolds of finite volume. Using the conformal covariance of the Dirac operator we prove:
\begin{theorem}\label{Hij1}
Let $(M^n, g)$ be a complete Riemannian $\Spinc$ manifold of finite volume and dimension $n>2$. Any eigenvalue $\lambda$ of the Dirac operator to which is attached an eigenspinor $\psi$ satisfies the Hijazi type inequality (\ref{ega}). Equality holds if and only if Equation (\ref{case}) holds.

\label{hjhjhj}
\end{theorem}
Now, the Hijazi type inequality (\ref{herr}) can be derived for complete Riemannian $\Spinc$ manifolds of finite volume and dimension $n>2$ and equality holds if and only if the eigenspinor associated with the first eigenvalue $\lambda_1$ is a Killing $\Spinc$ spinor.
This was also proved by N. Grosse in \cite{nadine} and \cite{nadine1} for complete spin manifolds of finite volume and dimension $n>2$. On complete manifolds, the Dirac operator is essentially self-adjoint and, in general, its spectrum consists of eigenvalues and the essential spectrum. For elements of the essential spectrum, we also extend to $\Spinc$ manifolds the Hijazi type inequality (\ref{herr}) obtained by N. Grosse in \cite{nadine1} on spin manifolds:
\begin{theorem}\label{n5sb}
Let $(M^n,g)$ be a complete Riemannian $\Spinc$ manifold of dimension $n\geq 5$ with finite volume. Furthermore, assume that $S -c_n \vert\Omega\vert$ is bounded from below. If $\lambda$ is in the essential spectrum of the Dirac operator $\sigma_{ess}(D)$, then $\lambda$ satisfies the Hijazi type  inequality (\ref{herr}). 
\end{theorem}
For the $2$-dimensional case, N. Grosse proved in \cite{nadine} that for any Riemannian spin surface of finite area, homeomorphic to $\R^2$ we have 
\begin{eqnarray}
 \lambda^+ \geq \frac{4\pi}{\aire},
\end{eqnarray}
Where $\lambda^+ = \inf_{\varphi \in C^\infty _c (M)} \frac{(D^2\varphi, \varphi)}{ (\varphi, \varphi)}$ (in the compact case, $\lambda^+$ coincides with the first eigenvalue of the square of the Dirac operator). Recently, in \cite{Bar}, C. B\"{a}r showed the same inequality for any connected $2$-dimensional Riemannian manifold of genus 0, with finite area and equipped with a spin structure which is {\it bounding at infinity}. A spin structure on $M$ is said to be {\it bounding at infinity} if $M$ can be embedded into $\SSS^2$ in such a way that the spin structure extends to the unique spin structure of $\SSS^2$.\\ \\
Studying the Energy-Momentum tensor on a compact Riemannian spin or $\Spinc$ manifolds has been done by many authors, since it is related to serval geometric situations. Indeed, on compact spin manifolds,  J. P. Bourguigon and P. Gauduchon \cite{bour-gaud, bgm} proved that the Energy-Momentum tensor appears naturally in the study of the variations of the spectrum of the Dirac operator. Hence, when deforming the Riemannian metric in the direction of this tensor, the eigenvalues of the Dirac operator are then critical. Using this, T. Friedrich and E.C. Kim \cite{fr-kim} obtained the Einstein-Dirac equation as the Euler-Lagrange equation of a certain functional. The author extends these last results to compact $\Spinc$ manifolds \cite{2ana}. Even it is not a computable geometric invariant, the Energy-Momentum tensor is, up to a constant, the second fundamental form of an isometric immersion  into a $\Spinc$ manifold carrying a parallel spinor \cite{morel, 2ana}. Moreover, in low dimensions the existence, on a spin or $\Spinc$ manifold $M$, of a spinor $\psi$ satisfying Equation (\ref{casee}) is, under some additional assumptions,  equivalent to the existence of a local immersion of $M$ into $\R^3$, $\SSS^3$, $\C P^2$, $\SSS^2\times \R$ or some others manifolds \cite{frfr, morel, rog-jul}.\\\\
{\bf Acknowledgment}\\\\
The author would like to thank Oussama Hijazi for his support and encouragements.
\section{Preliminaries}
In this section we briefly introduce basic notions concerning $\Spinc$ manifolds, the Dirac operator and its conformal covariance. Then we recall the refined Kato inequality which is crucial for the proof.\\ \\
{\bf The Dirac operator on $\Spinc$ manifolds:} Let $(M^n, g)$ be a connected  oriented Riemannian manifold of dimension $n\geqslant 2$ without boundary. Furthermore, let $SOM$ be the 
$SO_n$-principal bundle over $M$ of positively oriented orthonormal frames. A $\Spinc$ structure of $M$ is a $\Spin_n^c$-principal bundle $(\Spinc M,\pi,M)$
 and an $\cercle^1$-principal bundle $(\cercle^1 M ,\pi,M)$ together with a double covering given by  $\theta: \Spinc M \longrightarrow SOM\times_{M}\cercle^1 M$ such that
$\theta (ua) = \theta (u)\xi(a),$
for every $u \in \Spinc M$ and $a \in \Spin_n^c$, where $\xi$ is the $2$-fold covering of $\Spin_n ^c$ over $SO_n\times \cercle^1$. Let $\Sigma M := \Spinc M \times_{\rho_n} \Sigma_n$ be the associated spinor bundle where $\Sigma_n = \C^{2^{[\frac n2]}}$ and $\rho_n : \Spin_n^c
\longrightarrow  \End(\Sigma_{n})$ the complex spinor representation. A section of $\Sigma M$ will be called a spinor and the set of all spinors will be 
denoted by $\Gamma(\Sigma M)$ and those of compactly supported smooth spinors by $\Gamma_c(\Sigma M)$. The spinor bundle $\Sigma M$ is equipped with a natural Hermitian scalar product, denoted by $< . , .>$, satisfying
$$< X\cdot\psi, \varphi> = - <\psi, X\cdot\varphi> \ \ \text{for every}\ \ X\in \Gamma(TM)\ \ \text{and}\ \ \psi, \varphi\in \Gamma(\Sigma M),$$
where $X\cdot\psi$ denotes the Clifford multiplication of $X$ and $\psi$. With this Hermitian scalar product we define an $L^2$-scalar product
$$(\psi, \varphi) = \int_M <\psi, \varphi> \vo_g,$$
for any spinors $\psi$ and $\varphi$ in $\Gamma_c(\Sigma M)$. Additionally, given a connection 1-form $A$ on $\cercle^1 M$, $A: T(\cercle^1 M)\longrightarrow i\R$ and the connection 1-form 
$\omega^M$ on $SO M$ for the Levi-Civita connection $\nabla^M$, we consider the associated connection on the principal bundle $SO M\times_{M} \cercle^1 M$, and hence 
a covariant derivative $\nabla$ on $\Gamma(\Sigma M)$ \cite{5}.\\ 
The curvature of $A$ is an imaginary valued 2-form denoted by $F_A= dA$, i.e., $F_A = i\Omega$, where $\Omega$ is a real valued 2-form on $\cercle^1 M$. We know that $\Omega$ can be viewed as a real valued 2-form on $M$ $\cite{5}$. In this case $i\Omega$ is the curvature form of the associated line bundle $L$. It is the complex line bundle associated with the $\cercle^1$-principal bundle via the standard representation of the unit circle. For any spinor $\psi$ and any real $2$-form $\Omega$ we have \cite{8}:
\begin{eqnarray}
 <i\Omega\cdot\psi,\psi>\ \geqslant -\frac{c_n}{2} \vert \Omega \vert_g\vert \psi\vert^2,
\label{asese}
\end{eqnarray}
where $\vert \Omega \vert_g$ is the norm of  $\Omega$ given by $\vert \Omega \vert^2_g = \sum_{i<j} (\Omega_{ij})^2$ . Moreover,  equality holds in (\ref{asese}) if and only if 
\begin{eqnarray}
 \Omega\cdot\psi = i \frac {c_n}{2}\vert\Omega\vert_g \psi.
\label{asese=}
\end{eqnarray}
The Dirac operator is a first order elliptic operator locally given by  $$D =\sum_{i=1}^n e_i \cdot \nabla_{e_i}.$$ It is an elliptic and formally self-adjoint operator with respect to the $L^2$-scalar product, i.e., for all spinors $\psi$, $\varphi$ at least one of which is compactly supported on $M$ we have $(D\psi, \varphi) = (\psi, D\varphi)$. An important tool when examining the Dirac operator is the Schr\"{o}dinger-Lichnerowicz formula 
\begin{eqnarray}
 {D}^2 = {\nabla}^*\nabla + \frac 14 S\; \iid_{\Gamma (\Sigma M)}+ \frac{i}{2}\Omega\cdot,
\label{lich}
\end{eqnarray}
where ${\nabla}^*$  is the adjoint of $\nabla$ and  $\Omega\cdot$ is the extension of the Clifford multiplication to differential forms given by $(e_i ^* \wedge e_j ^*)\cdot\psi = e_i\cdot e_j \cdot\psi$. For the Friedrich connection $\nabla^f_X \psi = \nabla_X \psi + \frac{f}{n} X\cdot\psi$ where $f$ is real valued function one gets a Schr\"{o}dinger-Lichnerowicz type formula similar to the one obtained by Friedrich in \cite{friedrich}:
\begin{eqnarray}
 (D-f)^2\psi = \bigtriangleup^f \psi + (\frac{S}{4} +\frac{n-1}{n}f^2)\psi +\frac{i}{2} \Omega\cdot\psi- \frac{n-1}{n} (2f D\psi +\nabla f\cdot\psi),
\label{lfs}
\end{eqnarray}
where $ \bigtriangleup^f$ is the spinorial Laplacian associated with the connection $\nabla^f$.\\ 
A complex number $\lambda$ is an eigenvalue of $D$ if there exists a nonzero eigenspinor $\psi \in \Gamma(\Sigma M) \cap L^2 (\Sigma M)$ with $D\psi= \lambda\psi$. The set of all eigenvalues is denoted by $\sigma_p(D)$, the point spectrum. We know that if $M$ is closed, the Dirac operator has a pure point spectrum but on open manifolds, the spectrum might have a continuous part. In general the spectrum of the Dirac operator $\sigma (D)$ is composed of the point, the continuous and the residual spectrum. For complete manifolds, the residual spectrum is empty and $\sigma (D) \subset \R$. Thus, for complete manifolds, the spectrum can be divided into point and continuous spectrum. But often another decomposition of the spectrum is used: the one into discrete spectrum $\sigma_d(D)$ and essential spectrum $\sigma_{ess}(D)$.\\ 
 A complex number $\lambda$ lies in the essential spectrum of $D$ if there exists a sequence of smooth compactly supported spinors $\psi_i$ which are orthonormal with respect to the $L^2$-product and $$ \| (D-\lambda) \psi_i\|_{L^2} \longrightarrow 0.$$
The essential spectrum contains all eigenvalues of infinite multiplicity. In contrast, the discrete spectrum $\sigma_d(D) :=\sigma_p(D)\diagdown \sigma_{ess}(D)$ consists of all eigenvalues of finite multiplicity. The proof of the next property can be found in \cite{nadine}: On a $\Spinc$ complete Riemannian manifold, 0 is in the essential spectrum of $D-\lambda$ if and only if 0 is in the essential spectrum of $(D-\lambda)^2$ and in this case, there is a normalized sequence $\psi_i \in \Gamma_c(\Sigma M)$ such that $\psi_i$ converges $L^2$-weakly to $0$ with $\| (D-\lambda)\psi_i\|_{L^2} \longrightarrow 0$ and $\| (D-\lambda)^2\psi_i\|_{L^2} \longrightarrow 0$.\\ \\
{\bf Spinor bundles associated with conformally related metrics:} The conformal class of $g$ is the set of metrics $\overline g=e^{2u}g$, for a real  
function $u$ on $M$. At a given point $x$ of $M$, we consider a
$g$-orthonormal basis $\{e_1,\ldots,e_n\}\;$ of $T_xM$. The corresponding
$\overline g\,$-orthonormal basis is  denoted by 
$\{\overline{e}_1=e^{-u}e_1,\ldots,\overline{e}_n=e^{-u}e_n\}\;$. This correspondence extends to 
the $\Spinc$ level to give an isometry between the associated spinor
bundles. We put a ``\ $^{\overline{\;\;}}$\ '' above every object which is
naturally associated with the metric $\overline g.$ Then, for any 
spinor field $\psi $ and $\varphi $, one has
$<\overline{\psi} ,\overline{\varphi} >=<\psi ,\varphi >\,,$
where $<., .>$ denotes the natural Hermitian scalar products on 
$\Gamma (\Sigma M)$, and on $\Gamma (\Sigma \,\overline M)$.
The corresponding Dirac operators satisfy
\begin{eqnarray}
\overline D\ (\ e^{-\frac{n-1}{2}u}\;\overline \psi \ )=\ e^{-\frac{n+1}{2}u}\  \overline {D\psi}.
\label{Dconforme}
\end{eqnarray}
The norms of any real 2-form $\Omega$ with respect to $g$ and $\overline g$ are related by
\begin{eqnarray}
\vert\Omega\vert_{\overline g} \ = e^{-2u}\vert\Omega\vert_g.
\label{omegag}
\end{eqnarray}
O. Hijazi \cite {12} showed that on a spin manifold the Energy-Momentum tensor verifies 
\begin{eqnarray}
|\ell^{\overline{\varphi} }|_{\overline g}^2=e^{-2u}\,|\ell^\varphi |_g^2\,=e^{-2u}\,|\ell^\psi |_g^2,
\label{titi}
\end{eqnarray}
where $\varphi =e^{-\frac{(n-1)}{2}u}\psi$. We extend the result to a $\Spinc$ manifold and get the same relation.\\ \\
{\bf Refined Kato inequalities:} On a Riemannian manifold $(M, g)$, the Kato inequality states that away from the zeros of any section $\varphi$ of a Riemannian or Hermitian vector bundle $E$ endowed with a metric connection $\nabla$ we have,
\begin{eqnarray}
|d(|\varphi|)|\leq |\nabla\varphi|.
\label{kato}
\end{eqnarray}
This could be seen as follows $2|\varphi| |d(|\varphi|)|=|d(|\varphi|)^2|=2| <\nabla\varphi,\varphi>|\leq 2|\varphi||\nabla\varphi|$. In \cite{24}, refined Kato inequalities were obtained for sections in the kernel of first order elliptic differential operators $P$. They are of the form 
$ |d(|\varphi|)|\leq k_P|\nabla\varphi|,$ where $k_P$ is a constant depending on the operator $P$ and $0 < k_P < 1$. Without the assumption that $\varphi \in \ker P$, we get away from the zero set of $\varphi$
\begin{eqnarray}
|d|\varphi||\leq |P\varphi|+k_P|\nabla\varphi|.
\label{katofinal}
\end{eqnarray}
A proof of (\ref{katofinal}) can be found in \cite{24}, \cite{nadine}, \cite{branson} or \cite{nadine1}. In \cite{24} the constant $k_P$ is determined in terms of the conformal weights of the differential operator $P$. For the Dirac operator $D$ and for $D-\lambda$, where $\lambda\in \R$, we have $k_D = k_{D-\lambda} = \sqrt{\frac{n-1}{n}}$.
\section{Proof of the Hijazi type inequalities}
First, we follow the main idea of the proof of the original Hijazi inequality in the compact case (\cite{12}, \cite{13}), and its proof on spin noncompact case obtained by N. Grosse \cite{nadine1}. We choose the conformal factor with the help of an eigenspinor and we use cut-off functions near its zero-set and near infinity to obtain compactly supported test functions.
\begin{proof}[\bf {Proof of Theorem \ref{Hij1}}] Let $\psi\in C^\infty(M,S)\cap L^2(M,S)$ be a normalized eigenspinor, i.e., $D\psi=\lambda\psi$ and $\Vert \psi\Vert=1$. Its zero-set $\Upsilon$ is closed and lies in a closed countable union of smooth $(n-2)$-dimensional submanifolds which has locally finite $(n-2)$-dimensional Hausdorff measure \cite{Ba99}. We can assume without loss of generality that $\Upsilon$ is  itself a countable union of $(n-2)$-submanifolds described above. Fix a point $p\in M$. Since $M$ is complete, there exists a cut-off function $\eta_i: M\to [0,1]$ which is zero on $M\setminus B_{2i}(p)$ and equal $1$ on $B_{i}(p)$, where $B_{l}(p)$ is the ball of center $p$ and  radius $l$. In between, the function is chosen such that $|\nabla \eta_i|\leq \frac{4}{i}$ and $\eta_i\in C_c^\infty(M)$. While $\eta_i$ cuts off $\psi$ at infinity, we define another cut-off near the zeros of $\psi$. Let $\rho_{a,\epsilon}$ be the function

\[\rho_{a,\epsilon}(x)=\Bigg\{ \begin{array}{ll}  0  &\mathrm{for\ } r< a\epsilon\\ 
1-\delta\ln\frac{\epsilon}{r} &\mathrm{for\ } a \epsilon \leq r\leq\epsilon\\
1  &\mathrm{for\ } \epsilon < r \end{array}
 \]
where $r=d(x,\Upsilon)$ is the distance from $x$ to $\Upsilon$. The constant $ 0<a<1$ is chosen such that $\rho_{a,\epsilon} (a\epsilon)=0$, i.e., $a=e^{-\frac{1}{\delta}}$. Then $\rho_{a,\epsilon}$ is continuous, constant outside a compact set and Lipschitz.
Hence, for $\varphi\in \Gamma(\Sigma M)$ the spinor $\rho_{a,\epsilon}\varphi$ is an element in $ H_1^r(\Sigma M)$ for all $1\leq r\leq \infty$. Now, consider ${\Psi}:=\eta_i\rho_{a,\epsilon}\psi\in H_1^r(\Sigma M)$. These spinors are compactly supported on $M\setminus \Upsilon$. Furthermore, $\ov{g}=e^{2u}g=h^{\frac{4}{n-2}}g$ with $h=|\psi|^{\frac{n-2}{n-1}}$ is a metric on $M\setminus{\Upsilon}$. Setting $\ov{\Phi}:=e^{-\frac{n-1}{2}u}\ov{\Psi}$\ ($\varphi= e^{-\frac{n-1}{2}u}\psi$), Equations (\ref{asese}), (\ref{omegag}), (\ref{titi}) and the Schr\"{o}dinger-Lichnerowicz formula imply
\begin{eqnarray}
\|\overline\nabla^{\ell^{\overline{\Phi}}} \overline{{\Phi}}\|^2_{\overline g}
&=&  \|\overline D \ \overline {{ \Phi}}\|^2_{\overline g} \ - \frac14 \int_{M-\Upsilon}\overline S \vert\overline{{ \Phi}}\vert^2 \vo_{\overline g}\ -\int_{M-\Upsilon}\vert \ell^{\overline{{ \Phi}}}\vert^2\vert\overline{{ \Phi}}\vert^2 \vo_{\overline g} \nonumber\\ &&-\int_{M-\Upsilon} <\frac i2\Omega\ \overline \cdot\ \overline{{ \Phi}}, \overline{{ \Phi}}> \vo_{\overline g}\nonumber\\ & \leqslant &   
\|\overline D \ \overline {{ \Phi}}\|^2_{\overline g} - \frac 14\int_M ( \overline S e^{2u} -c_n \vert\Omega\vert_g) \vert{ \Psi}\vert^2  e^{-u} \vo_g - \int_M \vert \ell^{{ \Psi}}\vert^2  \vert{ \Psi}\vert^2 e^{-u}\vo_g\nonumber \\ &=& \|\overline D \ \overline {{ \Phi}}\|^2_{\overline g}  - \frac 14 \int_M (h^{-1} L^\Omega h) \vert{ \Psi}\vert^2 e^{-u} \vo_g - \int_M \vert \ell^{{ \Psi}}\vert^2  \vert{ \Psi}\vert^2 e^{-u}\vo_g,\nonumber
\end{eqnarray}
where $\nabla^{\ell^\varphi}_X\varphi$ is the spinor field defined in \cite{12} by $\nabla^{\ell^\varphi}_X\varphi:=\nabla_X\varphi+\ell^\varphi (X)\cdot\varphi$ and where we used  $|\ov{{ \Phi}}|^2\vo_{\ov{g}}=e^u|{ \Psi}|^2\vo_g$ and $\ov{S}e^{2u} - c_n\vert\Omega\vert_g =h^{-1}L ^\Omega h$ (see \cite{article1}). Using $\ov{D}\,\ov{\varphi} =\lambda e^{-u}\ov{\varphi}$ and $<\nabla (\eta_i\rho_{a,\epsilon})\overline \cdot\ \overline \varphi, \overline\varphi>\in  C^\infty (M, \i\R)$, we calculate 
\begin{eqnarray}
\|\overline D \ \overline {{ \Phi}}\|^2_{\overline g} =
\|\nabla (\eta_i\rho_{a,\epsilon}) \ \overline \cdot\  \overline \varphi \|^2_{\overline g} + \lambda^2 \int_M \eta_i^2 \rho_{a,\epsilon}^2 \ e^{-(n+2)u} \vert\varphi\vert^2 \vo_{g}.
\label{fafa}
\end{eqnarray}
Inserting (\ref{fafa}) and $\|\overline\nabla^{\ell^{\overline{{ \Phi}}}} \overline{{ \Phi}}\|^2_{\overline g} \geqslant 0$ in the above inequality, we get 
\begin{eqnarray*}
\Vert \nabla (\eta_i\rho_{a,\epsilon})\ \ov\cdot\ \ov{\varphi}\Vert^2_{\ov{g}}\geq \frac 14 \int_M (h^{-1} L^\Omega h) \vert{ \Psi}\vert^2 e^{-u} \vo_g  &+&\int_M \vert \ell^{{ \Psi}}\vert^2  \vert{ \Psi}\vert^2 e^{-u}\vo_g \\ &-& \lambda^2 \int_M \eta_i^2 \rho_{a,\epsilon}^2   \vert\varphi\vert^2 e^{-(n+2)u}\vo_{g}.
\end{eqnarray*}
Moreover, we have $\Vert \nabla (\eta_i\rho_{a,\epsilon})\ov\cdot\ \ov{\varphi}\Vert^2_{\ov{g}} =\int\limits_M |\nabla (\eta_i\rho_{a,\epsilon})\cdot\psi|^2e^{-u}\vo_{g}.$
Thus, with $e^u=|\psi|^\frac{2}{n-1}$ the above inequality reads
\begin{eqnarray*}
\int\limits_M |\nabla (\eta_i\rho_{a,\epsilon})|^2 |\psi|^{2\frac{n-2}{n-1}}\vo_g  &\geq& \frac{1}{4}\int\limits_M \eta_i\rho_{a,\epsilon} |\psi|^{\frac{n-2}{n-1}} L^\Omega (\eta_i\rho_{a,\epsilon}|\psi|^{\frac{n-2}{n-1}})\vo_g -\lambda^2\int\limits_M\eta_i^2\rho_{a,\epsilon}^2 |\psi|^{2\frac{n-2}{n-1}}\vo_g\\ &&-\frac{n-1}{n-2}\int\limits_M |\nabla (\eta_i\rho_{a,\epsilon})|^2 |\psi|^{2\frac{n-2}{n-1}}\vo_g + \int_M \vert \ell^\psi\vert^2 \vert \psi\vert^{2\frac{n-2}{n-1}} \eta_i^2\rho_{a,\epsilon}^2 \vo_g.
\end{eqnarray*} 
Hence, we obtain
\[
\frac{2n-3}{n-2}\!\int\limits_M\!\!|\nabla (\eta_i\rho_{a,\epsilon})|^2 |\psi|^{2\frac{n-2}{n-1}}\vo_g  \geq \left(\frac{\mu_1}{4} +\inf_M \vert \ell^\psi\vert^2 -\lambda^2\right)\!\int\limits_M\eta_i^2\rho_{a,\epsilon}^2 |\psi|^{2\frac{n-2}{n-1}}\vo_g,\]
where $\mu_1$ is the infimum of the spectrum of the perturbed conformal Laplacian.
With $|\eta_i \nabla \rho_{a,\epsilon} + \rho_{a,\epsilon} \nabla\eta_i|^2 \leq 2 \eta_i^2 |\nabla \rho_{a,\epsilon}|^2 + 2 \rho_{a,\epsilon}^2 |\nabla\eta_i|^2$ we have 
\[ k\!\int\limits_M\!(\eta_i^2|\nabla \rho_{a,\epsilon}|^2 +\rho_{a,\epsilon}^2|\nabla \eta_i|^2) |\psi|^{2\frac{n-2}{n-1}}\vo_g  \geq \left(\frac{\mu_1}{4} +\inf_M \vert \ell^\psi\vert^2-\lambda^2\right)\!\Vert \eta_i\rho_{a,\epsilon} |\psi|^{\frac{n-2}{n-1}}\Vert^2,\]
where $k=2\frac{2n-3}{n-2}$. 
Next, we examine the limits when $a$ goes to zero. Recall that $\Upsilon\cap \ov{B_{2i}(p)}$ is bounded, closed $(n-2)$-$C^\infty$-rectifiable and has still locally finite $(n-2)$-dimensional Hausdorff measure. For fixed $i$ we estimate
\[ \int\limits_M |\nabla \rho_{a,\epsilon}|^2\eta_i^2|\psi|^{2\frac{n-2}{n-1}}\vo_g\leq \sup_{B_{2i}(p)} |\psi|^{2\frac{n-2}{n-1}}\ \int\limits_{B_{2i}(p)} |\nabla \rho_{a,\epsilon}|^2\vo_g.\]
Furthermore, we set
$ \B_{\epsilon, p}:=\{x\in B_\epsilon\ |\ d(x,p)=d(x,{\Upsilon})\}$ with $B_\epsilon:=\{x\in M\ |\ d(x,{\Upsilon})\leq\epsilon\}$.
For $\epsilon$ sufficiently small each $ \B_{\epsilon, p}$ is star shaped. Moreover, there is an inclusion $\B_{\epsilon, p}\hookrightarrow B_\epsilon(0)\subset\R^{2}$ via the normal exponential map. Then we can calculate
\begin{align*}  \int\limits_{B_\epsilon\cap B_{2i}(p)}\hspace{-0.5cm} |\nabla \rho_{a,\epsilon}|^2\vo_g & \leq { \vol_{n-2}}
(\Upsilon\cap B_{2i}(p) ) \sup_{x\in\Upsilon\cap B_{2i}(p)} \int\limits_{\B_{\epsilon, x}\setminus \B_{a\epsilon, x}}
\hspace{-0.5cm} |\nabla \rho_{a,\epsilon}|^2\vo_{g'}\\
&\leq c \vol_{n-2}(\Upsilon\cap B_{2i}(p) ) \int\limits_{B_\epsilon(0)\setminus {B}_{a\epsilon}(0)}\hspace{-0.5cm} |\nabla \rho_{a,\epsilon}|^2\vo_{g_E}\\
& \leq c^\prime  \int\limits^\epsilon_{a\epsilon} \frac{\delta^2}{r}dr=-c^\prime \delta^2\ln a=c^\prime\delta\to 0 \quad \mathrm{for\ } a\to 0,
 \end{align*}
where $\vol_{n-2}$ denotes the $(n-2)$-dimensional volume 
and $g'=g_{|_{\B_{\epsilon, p}}}$. The positive constants $c$ and $ c^\prime$ arise from $\vol_{n-2}({\Upsilon}\cap B_{2i}(p))$ and the comparison of $\vo_{g'}$ with the volume element of the Euclidean metric. Furthermore, for any compact set $K \subset M$ and any positive function $f$ it holds $\rho_{a,\epsilon}^2 f \nearrow f$ and thus by the monotone convergence theorem, we obtain when $a \longrightarrow 0$, 
$$\int_K \rho_{a,\epsilon}^2 f \vo_g \longrightarrow \int_K f\vo_g.$$
When applied to the functions $\rho_{a,\epsilon}^2 |\nabla \eta_i|^2 |\psi|^{2\frac{n-2}{n-1}}$, with $K= B_{2i}(p)$ we get 
\[\int\limits_{B_{2i}(p)}\hspace{-0.2cm} \rho_{a,\epsilon}^2 |\nabla \eta_i|^2 |\psi|^{2\frac{n-2}{n-1}}\vo_g\to \int\limits_{B_{2i}(p)}\hspace{-0.2cm} |\nabla \eta_i|^2 |\psi|^{2\frac{n-2}{n-1}}\vo_g\]
 as $a\to 0$ and thus,
\[ k\int\limits_M |\nabla \eta_i|^2 |\psi|^{2\frac{n-2}{n-1}}\vo_g  \geq \left(\frac{\mu_1}{4} + \inf_M \vert \ell^\psi\vert^2-\lambda^2\right)\int\limits_M \eta_i^2 |\psi|^{2\frac{n-2}{n-1}}\vo_g.\\[0.2 cm] 
\]
Next we have to study the limit when $i\to \infty$: Since $M$ has finite volume and $\Vert\psi\Vert=1$, the H\"older inequality ensures that $\int\limits_M |\psi |^{2\frac{n-2}{n-1}}\vo_g$ is bounded. With $|\nabla \eta_i|\leq \frac{4}{i}$ we get the result. Equality is attained if and only if $\|\overline\nabla^{\ell^{\overline{{ \Phi}}}} \overline{{ \Phi}}\|^2_{\overline g}\longrightarrow 0$ for $i\to \infty$, $a\to 0$ and $\Omega\cdot\psi = i\frac {c_n}{2} \vert\Omega\vert_g \psi$. But we have
$$ 0\leftarrow \|\overline\nabla^{\ell^{\overline{{ \Phi}}}} \overline{{ \Phi}}\|^2_{\overline g}=
\Vert \eta_i\rho_{a,\epsilon}\overline\nabla^{\ell^{\overline{{ \Phi}}}} \overline\varphi +{\nabla}(\eta_i\rho_{a,\epsilon}) \ov\cdot\ \ \ov{\varphi} \Vert_{\ov{g}} \geq \Vert \eta_i\rho_{a,\epsilon}\overline\nabla^{\ell^{\overline{\varphi}}} \overline\varphi\Vert_{\ov{g}} -{\Vert {\nabla}(\eta_i\rho_{a,\epsilon}) \ov\cdot\ \ \ov{\varphi}\Vert_{\ov{g}}}.$$
Since $\Vert {\nabla}(\eta_i\rho_{a,\epsilon})\overline\cdot\ov{\varphi}\Vert_{\ov{g}}\to 0$, we conclude that $\overline\nabla^{\ell^{\overline{\varphi}}} \overline\varphi$ has to vanish on $M\setminus {\Upsilon}$.
\end{proof}
\begin{remark}
By the Cauchy-Schwarz inequality, we have 
\begin{eqnarray}
\vert \ell^\psi\vert^2 \geqslant \frac{(tr(\ell^\psi))^2}{n} =\frac{\lambda^2}{n},
\label{cs}
\end{eqnarray}
where $tr$ denotes the trace of $\ell^\psi$. Hence the Hijazi type inequality (\ref{herr}) can be derived. Equality is achieved if and only if the eigenspinor associated with the first eigenvalue $\lambda_1$ is a $\Spinc$ Killing spinor. In fact, if equality holds then $\lambda^2 = \frac{n}{4(n-1)} \mu_1 = \frac 14 \mu_1 + \vert \ell^\psi\vert^2$ and equality in (\ref{cs}) is satisfied. Hence it is easy to check that 
$$T^\psi (e_i, e_j) = 0 \ \text{for}\ \ i\neq j\ \ \text{and}\ \ T^\psi (e_i, e_i) = \pm \frac{\lambda}{n}.$$
Finally, $\ell^\psi(X) = \pm\frac{\lambda}{n}X$ and $\ell^{\overline \varphi} (X) = e^{-u} \ell^\psi(X)= \pm\frac{\lambda}{n} e^{-u}X$. By (\ref{case}) we get that $\overline \varphi$ is a generalized Killing $\Spinc$ spinor and hence a Killing $\Spinc$ spinor for $n\geqslant4$ (\cite[Theorem 1.1]{8}). The function $e^{-u}$ is then constant and $\psi$ is a Killing $\Spinc$ spinor. For $n=3$, we follow the same proof as in \cite{8}. First we suppose that $\lambda_1 \neq 0$, because if $\lambda_1 =0$, the result is trivial. We consider the Killing vector $\ov \xi$ defined by  $$ i\ov g (\ov \xi,X) = <X\ov\cdot\ \ov \varphi,\ov \varphi>_{\ov g}\ \ \ \ \text{for every}\ \ \ X \in \Gamma(TM).$$%
In \cite{8}, it is shown that $d\ov \xi  = 2\lambda_1 e^{-u} (*\ov \xi),$ $\nabla \vert \ov \xi\vert^2 = 0$ and 
$\ov\xi\  \ov\cdot\ \ov\varphi= i\vert \ov\xi\vert^2 \ov\varphi$, where $*$ is the Hodge operator defined on differential forms. Since $*\ov \xi (\ov \xi,.) =0$, the 2-form $\Omega$ can be written $\Omega=F \ov \xi +\ov\xi \wedge\alpha,$
where $\alpha$ is a real 1-form and $F$ a function. We have \cite{8}:
\begin{eqnarray}
\Omega(\ov\xi,.) &=& \vert\ov\xi\vert^2 \alpha(.) =-4\lambda_1 d(e^{-u}) (.),\\
\Omega\ \ov\cdot\ \ov\varphi &=& -i F\ \ov\varphi-i\vert\ov\xi\vert^2 \alpha\ \ov\cdot\ \ov\varphi.\nonumber
\label{jdjd}
\end{eqnarray}
But equality in (\ref{spinfr}) is achieved so $\Omega\ov\cdot\ \ov \varphi=i\frac{c_n}{2} \vert\Omega\vert_{\ov g} \ov\varphi,$ which implies that $\Omega\ov\cdot\ \ov \varphi$ is collinear to $\ov \varphi$ and hence $\alpha\ov \cdot\ \ov\varphi$ is collinear to $\ov \varphi$. Moreover, $d(e^{-u})(\ov\xi)=-\frac{1}{4\lambda_1}\Omega(\ov\xi,\ov\xi) =0$ so  $\alpha(\ov\xi)=0$. It is easy to check that  $<\alpha\ov\cdot\ \ov\varphi,\ov\varphi>_{\ov g}=0$ which gives $\alpha\ov\cdot\ \ov\varphi \perp \ov\varphi.$
Because of $\alpha\ov\cdot\ \ov\varphi \perp \ov\varphi$ and $\alpha\ov\cdot\ \ov\varphi$ is collinear to $\ov\varphi$, we have $\alpha\ov\cdot\ \ov\varphi =0$ and finally $\alpha=0.$ Using (\ref{jdjd}), we obtain $d(e^{-u})=0$, i.e., $e^{-u}$ is constant, hence $\ov\varphi$ is a Killing $\Spinc$ spinor and finally $\psi$ is also a $\Spinc$ Killing spinor.
\end{remark}
\begin{proof}[\bf Proof of Theorem \ref{ththth}] The proof of Theorem  \ref{ththth} is similar to Theorem \ref{Hij1}. It suffices to take $\ov g= g$, i.e., $e^u =1$. The Friedrich type inequality (\ref{spinfr}) is obtained from the Hijazi type inequality (\ref{herr}).
\end{proof}
Next, we want to prove Theorem \ref{n5sb} using the refined Kato inequality:
\begin{proof}[{\bf Proof of Theorem \ref{n5sb}}] We may assume $\vol(M,g)=1$.
If $\lambda$ is in the essential spectrum of $D$, then $0$ is in the essential spectrum of $D-\lambda$ and of $(D-\lambda)^2$. Thus, there is a sequence $\psi_i\in \Gamma_c(\Sigma M)$ such that $\Vert (D-\lambda)^2\psi_i\Vert\to 0$ and $\Vert(D-\lambda)\psi_i\Vert\to 0$ while $\Vert \psi_i\Vert=1$. We may assume that $|\psi_i|\in C_c^\infty(M)$. That can always be achieved by a small perturbation. Now let $\frac{1}{2}\leq \beta\leq 1$. Then $|\psi_i|^\beta\in H_1^2(M)$. First, we will show that the sequence $\Vert d (|\psi_i|^\beta)\Vert$ is bounded: By the Cauchy-Schwarz inequality we have
\begin{eqnarray*}
 \Big|\int\limits_{|\psi_i|\ne 0} |\psi_i|^{2\beta-2}<(D-\lambda)^2\psi_i,\psi_i>\vo_g\Big| &\leq& \Vert |\psi_i|^{2\beta-1}\Vert_{\{|\psi_i|\ne 0\}} \Vert (D-\lambda)^2\psi_i\Vert \\ &\leq& \Vert\psi_i\Vert^{2\beta-1} \Vert (D-\lambda)^2\psi_i\Vert= \Vert (D-\lambda)^2\psi_i\Vert.
\end{eqnarray*}
Using (\ref{asese}) and the Schr\"{o}dinger-Lichnerowicz type formula (\ref{lfs}), we obtain
\begin{eqnarray*}
\Vert (D-\lambda)^2\psi_i\Vert
&\geq &\int\limits_{|\psi_i|\ne 0}\hspace{-0.3cm} |\psi_i|^{2\beta-2}|\nabla^\lambda\psi_i|^2\vo_g+
\!2(\beta-1)\!\int\limits_{|\psi_i|\ne 0}\hspace{-0.3cm} |\psi_i|^{2\beta-3}\!<d|\psi_i|\!\cdot\!\psi_i,\nabla^\lambda\psi_i\!>\!\vo_g\\
&&
+\!\int\!\left(\frac{S}{4}- \frac{c_n}{4} \vert\Omega\vert-\frac{n-1}{n}\lambda^2\right) |\psi_i|^{2\beta}\vo_g
\\ && -2\frac{n-1}{n}\lambda\Vert |\psi_i|^{2\beta-1}\Vert_{\{|\psi_i|\ne 0\}} \Vert(D-\lambda)\psi_i\Vert.
\end{eqnarray*}
The Cauchy-Schwarz inequality and the refined Kato inequality (\ref{kato}) for the connection $\nabla^\lambda$ imply
\begin{eqnarray*}
 && \int\limits_{|\psi_i|\ne 0}\hspace{-0.3cm} |\psi_i|^{2\beta-2}|\nabla^\lambda\psi_i|^2\vo_g + \!2(\beta-1)\!\int\limits_{|\psi_i|\ne 0}\hspace{-0.3cm} |\psi_i|^{2\beta-3}\!<d|\psi_i|\!\cdot\!\psi_i,\nabla^\lambda\psi_i\!>\!\vo_g 
\\ &\geq& (2\beta-1) \int\limits_{|\psi_i|\ne 0}\hspace{-0.3cm} |\psi_i|^{2\beta-2}|d(|\psi_i|)|^2\vo_g =(2\beta-1)\frac{1}{\beta^2}\int\limits_{|\psi_i|\ne 0} |d(|\psi_i|^\beta)|^2\vo_g.
\end{eqnarray*}
Hence, we have
\begin{eqnarray*}
\Vert (D-\lambda)^2\psi_i\Vert
&\geq &(2\beta-1)\frac{1}{\beta^2}\int\limits_{|\psi_i|\ne 0} |d(|\psi_i|^\beta)|^2\vo_g+ \int\limits\left(\frac{S}{4}-\frac{c_n}{4}\vert\Omega\vert-\frac{n-1}{n}\lambda^2\right) |\psi_i|^{2\beta}\vo_g\\
&&-2\frac{n-1}{n}\lambda \Vert(D-\lambda)\psi_i\Vert.
\end{eqnarray*}
Since $S-c_n \vert\Omega\vert$ is bounded from below, $\int\limits (S-c_n \vert\Omega\vert)|\psi_i|^{2\beta}\vo_g\geq \inf (S-c_n \vert\Omega\vert)\ \Vert \psi_i \Vert^{2\beta}_{2\beta}\geq \min\{\inf (S-c_n \vert\Omega\vert), 0\} $ is also bounded. Thus, with $\Vert (D-\lambda)\psi_i\Vert\to 0$ we see that $\Vert d|\psi_i|^\beta\Vert$ is also bounded. Next we fix $\alpha=\frac{n-2}{n-1}$ and obtain
\begin{eqnarray*}
\frac{\mu_1}{4}-\frac{n-1}{n}\lambda^2 &\leq& \left(\frac{\mu_1}{4}-\frac{n-1}{n}\lambda^2\right)\Vert |\psi_i|^{\alpha}\Vert^2\\
& \leq&  \frac{1}{4}\int\limits |\psi_i|^\alpha L^\Omega |\psi_i|^\alpha\vo_g-\frac{n-1}{n}\lambda^2\Vert|\psi_i|^\alpha\Vert^2\\
&=&\!\int\!\! |\psi_i|^{2\frac{n-2}{n-1}-2}\Big[(\frac{n}{n-1} |d(|\psi_i|)|^2+\frac{1}{2}d^*d(|\psi_i|^2)
\\&&+\left(\frac{S}{4}-\frac{c_n}{4}\vert\Omega\vert-\frac{n-1}{n}\lambda^2\right)\!|\psi_i|^2\Big]\vo_g,
\end{eqnarray*}
where we used the definition of $\mu_1$ as the infimum of the spectrum of $L^\Omega$ and  $|\psi_i|^\alpha d^*d(|\psi_i|^\alpha) = \frac{\alpha}{2}|\psi_i|^{2\alpha-2}d^*d(|\psi_i|^2)-\alpha(\alpha-2)|\psi_i|^{2\alpha-2}|d(|\psi_i|)|^2.$
Next, using the following:
$$\frac{1}{2}d^*d<\psi_i,\psi_i> \ \leqslant\  <D^2\psi_i,\psi_i>-\frac 14(S -{c_n}\vert\Omega\vert)|\psi_i|^2 -|\nabla\psi_i|^2,$$
$$|\nabla^\lambda\psi_i|^2 =|\nabla\psi_i|^2-2\frac{\lambda}{n}\Re <(D-\lambda)\psi_i,\psi_i>-\frac{\lambda^2}{n}|\psi_i|^2,$$
we have
\begin{eqnarray*}
\frac{\mu_1}{4}-\frac{n-1}{n}\lambda^2 
&\leq& \int\limits |\psi_i|^{2\frac{n-2}{n-1}-2}\left(\frac{n}{n-1} |d(|\psi_i|)|^2-|\nabla^\lambda\psi_i|^2\right)\vo_g\\
&&+\int\limits |\psi_i|^{2\frac{n-2}{n-1}-2}<(D-\lambda)^2\psi_i,\psi_i>\vo_g\\
&&+ \int\limits 2\left(1-\frac{1}{n}\right)\lambda |\psi_i|^{2\frac{n-2}{n-1}-2}\Re<(D-\lambda)\psi_i,\psi_i>\vo_g.
\end{eqnarray*}
The limit of the last two summands vanish since
$$\left|\int\limits |\psi_i|^{2\frac{n-2}{n-1}-2}<(D-\lambda)^2\psi_i,\psi_i>\vo_g\right|\leq \Vert (D-\lambda)^2\psi_i\Vert\ \Vert\, |\psi_i|^{\frac{n-3}{n-1}}\Vert\to 0,$$
$$\Big|\int\limits |\psi_i|^{2\frac{n-2}{n-1}-2}\Re<(D-\lambda)\psi_i,\psi_i>\vo_g\Big|\leq \Vert (D-\lambda)\psi_i\Vert\ \Vert\, |\psi_i|^{\frac{n-3}{n-1}}\Vert\to 0.$$
For the other summand we use the Kato-type inequality (\ref{katofinal}), $$|d(|\psi|)|\leq |(D-\lambda)\psi|+k |\nabla^\lambda\psi|,$$ which holds outside the zero set of $\psi$ and where $k=\sqrt{\frac{n-1}{n}}$. Thus, for $n\geq 5$ we can estimate 
\begin{eqnarray*}
&&\int\limits |\psi_i|^{2\alpha-2}\left(\frac{n}{n-1} |d(|\psi_i|)|^2 -|\nabla^\lambda\psi_i|^2\right)\vo_g\\
&=& \int\limits |\psi_i|^{2\alpha-2}\left(k^{-1} |d(|\psi_i|)|-|\nabla^\lambda\psi_i|\right)\left(k^{-1} |d(|\psi_i|)|+|\nabla^\lambda\psi_i|\right)\vo_g\\
&\leq& k^{-1}\hspace{-0.2cm}\int\limits_{\{|d(|\psi_i|)|\geq k|\nabla^\lambda\psi_i|\}}\hspace{-0.2cm} |\psi_i|^{2\alpha-2}|(D-\lambda)\psi_i|\big(k^{-1} |d(|\psi_i|)|+|\nabla^\lambda\psi_i|\big)\vo_g\\
&\leq& 2k^{-2}\hspace{-0.2cm}\int\limits_{\{|d(|\psi_i|)|\geq k|\nabla^\lambda\psi_i|\}} |\psi_i|^{2\alpha-2}|(D-\lambda)\psi_i| |d(|\psi_i|)|\vo_g\\
&\leq& 2k^{-2}\frac{n-1}{n-3} \Vert (D-\lambda)\psi_i\Vert\ \Vert d(|\psi_i|^{\frac{n-3}{n-1}})\Vert.
\end{eqnarray*}
For $n\geq 5$ we have $1\geq\frac{n-3}{n-1}\geq \frac{1}{2}$ and, thus, $\Vert d|\psi_i|^{\frac{n-3}{n-1}}\Vert$ is bounded. Together with $\Vert (D-\lambda)\psi_i\Vert\to 0$ we obtain the following: For all $\epsilon>0$ there is an $i_0$ such that for all $i\geq i_0$ we have
\[\int\limits |\psi_i|^{2\frac{n-2}{n-1}-2}\left(\frac{n}{n-1} |d|\psi_i||^2-|\nabla^\lambda\psi_i|^2\right)\vo_g\leq \epsilon.\]
Hence, we have $\frac{\mu_1}{4}\leq \frac{n-1}{n}\lambda^2$.
\end{proof}

\end{document}